\documentclass[11pt,reqno]{amsart}
\usepackage{amsmath}
\usepackage{amssymb}
\usepackage{amsfonts}
\usepackage{euscript}
\usepackage[latin1]{inputenc}
\usepackage[T1]{fontenc}

\numberwithin{equation}{section}
\makeatletter\@addtoreset{equation}{section}

\DeclareMathSymbol{\subsetneqq}{\mathbin}{AMSb}{36}
\newcommand{\set}[1]{{\left\{{#1}\right\}}}

\newcommand{\C}{\mathbb C}
\newcommand{\R}{\mathbb R}
\newcommand{\Z}{\mathbb Z}

\newtheorem {Theorem}{Main Theorem}[section]

\newtheorem {Proposition}[Theorem]{Proposition}
\newtheorem {Remark}[Theorem]{Remark}

\newtheorem {Claim}{Claim}

\begin{document}

\centerline{}

\title[On $L^2$-eigenfunctions of Twisted Laplacian on
 curved surfaces]{\texttt{On $L^2$-eigenfunctions of Twisted Laplacian on
 curved surfaces and suggested orthogonal polynomials}}

\author{A. Ghanmi}
\address{Department of Mathematics,  Faculty of  Sciences, P.O. Box 1014\newline
   Mohammed V University,  Agdal,  10 000 Rabat - Morocco  }
  \email{allalghanmi@gmail.com}

\maketitle

\begin{abstract} We show in a unified manner that the factorization method describes
 completely the $L^2$-eigenspaces associated to the discrete
 part of the spectrum of the twisted Laplacian on constant
curvature Riemann surfaces. Subclasses of two variable orthogonal polynomials are then derived and
arise by successive derivations of elementary complex valued
functions depending on the geometry of the surface.
\end{abstract}

\section{ Introduction and preliminaries.}
Let $M_\kappa$ be a given simply connected Riemann surface of
constant scalar curvature $\kappa\in\R$ (SCRS for short).
Precisely $M_\kappa$ is the disc of radius $1/\sqrt{-\kappa}$ for
$\kappa <0$, the Euclidean plane for $\kappa=0$ and the sphere in
$\R^3$ of radius $1/\sqrt{\kappa}$ identified with the extended
complex plane $\C\cup\{\infty\}$ for $\kappa
>0$. The corresponding Bergman-K\"ahler geometry is the one described by
the Hermitian metric $ ds^2_{\kappa} := (1+\kappa z \bar z)^{-2}dz
\otimes d\bar z,$
 $z=x+iy \in \C,$ whose  associated  volume
measure is $$ d\mu_{\kappa} = (1+\kappa z \bar z)^{-2}d\lambda,$$
where $d\lambda$ denotes the usual Lebesgue measure on $M_{\kappa}$.
Also let $\theta_\kappa$ be the differential one form
$\theta_\kappa:= (1+\kappa z \bar z)^{-1}(\bar zdz - zd\bar z).$
Associated to $M_\kappa$, $ds^2_{\kappa}$, $d\mu_{\kappa}$ and
$\theta_\kappa$, we consider the twisted Laplacian
$\mathfrak{L}_\kappa^\nu $, $\nu>0$, realized as magnetic
Schr\"odinger operator through $$\mathfrak{L}_\kappa^\nu = (d+i\nu
\theta_\kappa)^{*}(d+i\nu\theta_\kappa), $$ and acting on the
Hilbert space $\mathcal{H}_\kappa :=L^2(M_\kappa ;d\mu_{\kappa})$.
It is an elliptic self-adjoint second order differential operator
describing a single non relativistic spineless particle constrained
to move on the two-dimensional analytic surface $M_\kappa$ in the
presence of the
 external constant magnetic field  ${\mathcal{B}}= \nu d\theta_\kappa$.
The explicit expression of $\mathfrak{L}_\kappa^\nu$ in the $z$-complex
variable is given (up to a multiplicative constant)  by
\begin{eqnarray}
\mathfrak{L}_\kappa^\nu= -\{(1+\kappa |z|^2)^2\frac {\partial^2} {\partial z\partial \bar z} +\nu (1+\kappa
|z|^2)(z\frac {\partial}{\partial z}- {\bar z}\frac{\partial}{\partial{\bar z}}) -\nu^2 |z|^2\}.
\label{LH}
\end{eqnarray}
Its concrete spectral theory is well known in the literature, see
for example \cite{LL,AP} for $\kappa=-1,0,+1$ and \cite{FV,GhInJMP}
for arbitrary $\kappa$. In particular, the discrete part of the
spectrum of $\mathfrak{L}_\kappa^\nu$  acting on
$\mathcal{H}_\kappa$ is given by the eigenvalues $$
E_{\kappa,m}^\nu:=\nu(2m+1) + m(m+1)  \kappa$$ with $m$ is a
positive integer such that $0\leq m< ({2\nu +
\kappa})/({|\kappa|-\kappa})$ and the conditions $2\nu +\kappa
>0$ for $\kappa \leq 0$ and $2\nu/\kappa$ is integer for $\kappa>0$.
Moreover, we have the following (see \cite{GhJPA,PZ,GhInJMP} for example):

\begin{Proposition} \label{PropLH}  An orthogonal basis of the $L^2$-eigenspace
  $$A^{2,\nu}_m(M_\kappa):=\{\Phi \in \mathcal{H}_\kappa;
  ~~ \mathfrak{L}_\kappa^\nu  \Phi =E_{\kappa,m}^\nu \Phi  \}$$
  is given in terms of the real Jacobi Polynomials
  $\mathrm{P}^{(\alpha,\beta)}_{l}(x)$ by
   $$ (1+\kappa |z|^2)^{-(\frac \nu\kappa +m)} z^p\bar z^q
  \mathrm{P}^{(p+q,-2(\frac \nu\kappa +m)-1)}_{m-q}(1+2\kappa |z|^2)$$ with $q\leq m$ and the convention that $pq=0$.
\end{Proposition}

In the other hand, we know that the factorization algebraic method \cite{Sch40/41,Infeld51,FV} allows one
to construct $L^2$-eigenfunctions of second order differential operator like
$\mathfrak{L}_\kappa^\nu$. Our goal here is to discuss the converse. Precisely, we have to show
that such method describes completely the $L^2$-eigenspaces
associated to the discrete part of the twisted Laplacians $\mathfrak{L}_\kappa^\nu$ (following the terminology of \cite{AP,Ricci}). This will be done in a unified manner taking into account the curvature of the considered SCRS. In
particular, we recover the case $\kappa=+1$ discussed in \cite{FV}. The suggested subclass of two variable orthogonal
polynomials $P_{m,n}^{\nu; \kappa}(z, \bar z)$, i.e., such that the
functions $$(1+\kappa |z|^2)^{-(\frac \nu \kappa +m)} P_{m,n}^{\nu;
\kappa}(z, \bar z)$$ are $L^2$-eigenfunctions of
$\mathfrak{L}_\kappa^\nu$ with $E_\kappa^\nu(m):=
\nu(2m+1)+m(m+1)\kappa$ as corresponding $L^2$-eigenvalue, are
derived and satisfy the following Rodrigues type formula
\begin{eqnarray*} P_{m,n}^{\nu; \kappa}(z, \bar z)= (-1)^m
(1+ \kappa |z|^2)^{2\frac{\nu}{\kappa}+m} [(1+\kappa |z|^2)^{2}
\frac{\partial}{\partial z}]^{m}(1+ \kappa
|z|^2)^{-2(\frac{\nu}{\kappa}+m)}z^n, \label{KappaPoly0}\end{eqnarray*}
 where the intertwining  invariant
operator $[(1+ \kappa |z|^2)^{2}\frac{\partial}{\partial z}]^m$
depends only on the geometry of the considered SCRS  and not on the
magnetic field.

\section{ Main result.} We begin by recalling briefly the factorization method for curved surfaces. Indeed, one considers the first order
differential operator $\nabla_\alpha $ and its formal
adjoint $ \nabla_\alpha^{*}$ given respectively by $$ \nabla_\alpha = -(1+\kappa
|z|^2)\frac{\partial}{\partial z} + \alpha \bar z, \qquad\qquad  \nabla_\alpha^{*}  = (1+\kappa
|z|^2)\frac{\partial}{\partial \bar z} + (\alpha -\kappa) z.$$ Thus,
 by direct computation one gets the  algebraic relationships:
 \begin{align*}
\nabla_{\nu+\kappa}  \nabla_{\nu+\kappa}^{*}  = \mathfrak{L}_\kappa^\nu -\nu \quad \mbox{and} \quad
\nabla_{\nu+\kappa}^{*}  \nabla_{\nu+\kappa}  = \mathfrak{L}_\kappa^{\nu +\kappa} +(\nu+\kappa) ,
 \end{align*}
   which gives rise to the following one
$$ \mathfrak{L}_\kappa^\nu \nabla_{\nu+\kappa}=  \nabla_{\nu+\kappa}
\mathfrak{L}_\kappa^{\nu +\kappa} + (2\nu+\kappa)\nabla_{\nu+\kappa}. $$
Hence, for $\Psi_{0}$  being a nonzero $L^2$-eigenfuncion associated to the  lowest Landau level
 $E_{\kappa,0}^{\nu+ m \kappa}=\nu+ m \kappa$  of $\mathfrak{L}_\kappa^{\nu+ m \kappa}$,
one can show that $\nabla_{\nu + m \kappa} \Psi_{0}$ is also an
$L^2$-eigenfunction but this time for $\mathfrak{L}_\kappa^{\nu +
(m-1)\kappa}$ with
 $E_{\kappa,1}^{\nu + (m-1)\kappa}$ as corresponding eigenvalue.
 Doing so, it follows that
\begin{eqnarray} \nabla_{\nu+\kappa} \circ
\nabla_{\nu+2\kappa}\circ \cdots \circ \nabla_{\nu + m \kappa}
\Psi_{0} \in A^{2,\nu}_m(M_\kappa)
\label{ClosedEigenfunction}\end{eqnarray} i.e., it is an
$L^2$-eigenfunction of $\mathfrak{L}_\kappa^\nu $ associated to the
  eigenvalue $\nu(2m+1)+m(m+1)\kappa =:E_{\kappa,m}^\nu.$
 Conversely,
  we  show that every  $L^2$-eigenfunction $\Phi\in
A^{2,\nu}_m(M_\kappa)$
 is obtained by (\ref{ClosedEigenfunction}), that is each $\Phi\in
A^{2,\nu}_m(M_\kappa)$ arises by successive derivations of
 an elementary complex valued function. More precisely, we have

\begin{Theorem} \label{MainTheorem}  Fix $\nu>0$ such
that $2\nu + \kappa>0$ for $\kappa\leq 0$ and $2\nu/\kappa\in \Z^+$
for $\kappa>0$, and let $m$ be a fixed positive integer satisfying
$0\leq m< ({2\nu + \kappa})/({|\kappa|-\kappa})$. Then, the
$L^2$-eigenfunctions \begin{eqnarray} \label{ClosedEigenfunction2} \nabla_{\nu+\kappa} \circ
\nabla_{\nu+2\kappa}\circ \cdots \circ \nabla_{\nu + m
\kappa}[(1+\kappa |z|^2)^{-(\frac{\nu}{\kappa}+m)}z ^{n}], \quad
n=0,1,2, \cdots ,
\end{eqnarray}
 constitute an orthogonal basis of the $L^2$-eigenspace
$A^{2,\nu}_m(M_\kappa)$.
\end{Theorem}

\begin{Remark}  The above result says that the factorization
method determines completely all $L^2$-eigenfunctions of
$\mathfrak{L}_\kappa^\nu$, i.e., solutions of the eigenvalue problem
$\mathfrak{L}_\kappa^\nu \Psi =E_{\kappa,m}^\nu\Psi$ associated to
the discrete part of the $L^2$-spectrum. The planar case
($\kappa=0$) is classic. For $\kappa=+1$ (i.e., for the sphere
$S^2\cong\C\cup\{\infty\}$ equipped with the Fubini-Study metric on
the chart $\C$) the result has been established by Ferapontov and
Veselov \cite[Theorem 3]{FV}. While when $\kappa=-1$  the result we
obtain  can be considered as its analogue for the non compact
hyperbolic unit disc.
 \end{Remark}

The proof of the above theorem relies essentially on Proposition
\ref{PropLH} and the following result giving closed explicit
expressions of \eqref{ClosedEigenfunction} or also \eqref{ClosedEigenfunction2}. Namely, we have

\begin{Proposition}  \label{GhanmiProp}  Fix $\nu$
and $m$ as in the theorem above
 and define $P_{m,n}^{\nu; \kappa}(z, \bar z)$
by
$$P_{m,n}^{\nu; \kappa}(z, \bar z):=(1+\kappa |z|^2)^{\frac{\nu}{\kappa}+m}\nabla_{\nu+\kappa} \circ
\nabla_{\nu+2\kappa}\circ \cdots \circ \nabla_{\nu + m \kappa}[(1+\kappa |z|^2)^{-(\frac{\nu}{\kappa}+m)}z ^{n}]
$$  denote $m\wedge n:= Min(m,n)$ and set
\[C_{\kappa,\nu}^{m,n}= {(-1)^{m+n}} \frac{\Gamma(2(\frac{\nu}{\kappa}+m)-(m+n)+1)}
{\kappa^{n}\Gamma(2(\frac{\nu}{\kappa}+m)-m+1)}
.\]
 Then, we have
\begin{align} P_{m,n}^{\nu; \kappa}(z,\bar z) 
&=      (-1)^m(1+ \kappa |z|^2)^{2\frac{\nu}{\kappa}+m} [(1+\kappa |z|^2)^{2} \frac{\partial}{\partial z}]^{m}\Big((1+
\kappa |z|^2)^{-2(\frac{\nu}{\kappa}+m)}z^n \Big)\label{RTF1}
\\ &=    C_{\kappa,\nu}^{m,n}\cdot(1+\kappa |z|^2)^{2(\frac{\nu}{\kappa}+m)+1} \frac{\partial^{m+n}}{\partial z^m
\partial \bar z ^n}(1+ \kappa |z|^2)^{-2(\frac{\nu}{\kappa}+m)+m+n-1} \label{RTF2}
\\ &=   (-1)^{m}  (m \wedge n)!
|z|^{|m-n|}e^{i[(n-m)\arg z]} \mathrm{P}^{(|m-n|,-2(\frac{\nu}{\kappa}+m)-1)}_{m \wedge
n}(1+2\kappa |z|^2). \quad\label{RTF3}
 \end{align}
\end{Proposition}

\begin{proof}[Sketched proof  of Proposition \ref{GhanmiProp}]
 The identity (\ref{RTF1}) holds by observing
that the first order differential operator $\nabla_\alpha=-(1+\kappa
|z|^2)\frac{\partial}{\partial z} +\alpha \bar z$ can be rewritten
 as
\[ \nabla_{\alpha} f=-(1+\kappa |z|^2)^{\frac \alpha
\kappa+1}\frac{\partial}{\partial z}
 [(1+\kappa |z|^2)^{-\frac \alpha \kappa} f ]  \]
for every smooth function $f$ on $M_\kappa$. Thus, we have
\begin{align*}
\nabla_{\nu+\kappa} \circ \nabla_{\nu+2\kappa}\circ \cdots \circ
\nabla_{\nu +m\kappa}f
=(-1)^m(1+\kappa z\bar z)^{\frac \nu\kappa }  [(1+\kappa z\bar
z)^{2}\frac{\partial}{\partial z} ]^m ((1+\kappa z\bar z)^{-(\frac
\nu\kappa +m)} f).
\end{align*}
This yields (\ref{RTF1}) when specifying  $f(z) = (1+\kappa z\bar
z)^{-(\frac{\nu}{\kappa}+m)}z ^{n}$.

 The identity (\ref{RTF2}) is deduced from (\ref{RTF1}) using
\begin{Claim} \label{lemVeryClosed}
 For every fixed positive integer $m$ and every
smooth complex valued function $f$ on $M_\kappa$, we have
$$[(1+\kappa |z|^2)^2\frac{\partial}{\partial z}]^m f=(1+\kappa
|z|^2)^{m+1} \frac{\partial^m}{\partial z^m} ( (1+\kappa
|z|^2)^{m-1}f) $$
\end{Claim}
\noindent combined with the fact that
$$
z^n (1+\kappa
 |z|^2)^\alpha =\frac{1}{(\alpha+1)_n
 \kappa^n}
 \frac{\partial^{
 n}}{ \partial \bar z^n}
 ((1+\kappa
 |z|^2)^{\alpha+n}) ,$$ keeping in mind that $(-a)_n=(a -n+1)_n$. Here $(a)_n=a(a+1) \cdots (a+n-1)$.

The proof of (\ref{RTF3}) can be handled by induction together with
the use of the following result satisfied by the  real Jacobi Polynomials
  $\mathrm{P}^{(a,b)}_{j}(x)$:

\begin{Claim} \label{handbook} For arbitrary real number $a,b$, we
have $$(x^2-1)\frac{d}{d x}P^{(a,b)}_j(x) +
[(a-b)+(a+b)x]P^{(a,b)}_j(x)= 2(j+1)P^{(a-1,b-1)}_{j+1}(x).
$$
\end{Claim}
\end{proof}

Below, we give the proofs of Claims \ref{lemVeryClosed} and \ref{handbook}.

\begin{proof}[Proof of Claim \ref{lemVeryClosed}] This is proved by induction, where $m=0$ and $m=1$
is obvious. Next, set $ { \not \hspace*{-.1cm}D_\kappa}f=
h^2\frac{\partial}{\partial z} f = (1+\kappa z\bar
z)^2\frac{\partial}{\partial z}f,$  let ${ \not
\hspace*{-.1cm}D_\kappa^m}$ stand for $$ { \not
\hspace*{-.1cm}D_\kappa^m} = \underbrace{{ \not
\hspace*{-.1cm}D_\kappa}\circ { \not \hspace*{-.1cm}D_\kappa}\circ
\cdots \circ { \not \hspace*{-.1cm}D_\kappa}}_{\mbox{m-times}}$$ and
assume that $$ h^{m+1} \frac{\partial^m}{\partial z^m} ( h^{m-1}f)=
(h^2\frac{\partial}{\partial z})^m f ={ \not
\hspace*{-.1cm}D_\kappa^m} f $$ is satisfied for a given positive
integer $m$.  Hence, using direct computation combined with the fact
that $ { \not \hspace*{-.1cm}D_\kappa^m}((\frac{\partial h}{\partial
z}) f) = (\frac{\partial h}{\partial z}) { \not
\hspace*{-.1cm}D_\kappa^m} f $ for $\frac{\partial h}{\partial z}$
being an holomorphic function, we get \begin{eqnarray} h^{m+2}
\frac{\partial^{m+1}}{\partial z^{m+1}} (
h^{m}f)=mh\frac{\partial}{\partial z}(h{ \not
\hspace*{-.1cm}D_\kappa^m} f )+ h { \not
\hspace*{-.1cm}D_\kappa^m}(h^{-1}{ \not \hspace*{-.1cm}D_\kappa}f).
\label{RecF1}
\end{eqnarray} Next, note that we have \begin{eqnarray}   \qquad \quad h {
\not \hspace*{-.1cm}D_\kappa^m}(h^{-1}{ \not
\hspace*{-.1cm}D_\kappa}f) = h { \not
\hspace*{-.1cm}D_\kappa^{m-1}}(h^{-1}{ \not
\hspace*{-.1cm}D_\kappa}^2f) - h\frac{\partial h}{\partial z}{ \not
\hspace*{-.1cm}D_\kappa^{m-1}}({ \not \hspace*{-.1cm}D_\kappa}f).
\label{3Terms} \end{eqnarray} Repeated application of (\ref{3Terms})
gives
$$h { \not \hspace*{-.1cm}D_\kappa^m}(h^{-1}{ \not \hspace*{-.1cm}D_\kappa}f) = h
{ \not \hspace*{-.1cm}D_\kappa^{m-j}}(h^{-1}{ \not
\hspace*{-.1cm}D_\kappa^{j+1}}f) - j h\frac{\partial h}{\partial z}{
\not \hspace*{-.1cm}D_\kappa^{m-1}}({ \not
\hspace*{-.1cm}D_\kappa}f)$$ for every given positive integer $j$
such that $0\leq j\leq m$. In particular, for $j=m$ it follows
\begin{eqnarray}  \qquad \qquad  m h\frac{\partial h}{\partial z}{ \not
\hspace*{-.1cm}D_\kappa}^{m}(f) + h { \not
\hspace*{-.1cm}D_\kappa^m}(h^{-1}{ \not \hspace*{-.1cm}D_\kappa}f)={
\not \hspace*{-.1cm}D_\kappa^{m+1}}.\label{RecF2}
\end{eqnarray}
Finally, by combining (\ref{RecF1}) and (\ref{RecF2}), we get the
desired result of Claim \ref{lemVeryClosed}.
\end{proof}

\begin{proof}[Proof of Claim \ref{handbook}]
 The assertion of Claim \ref{handbook} is an immediate consequence of the facts that the
  classical Jacobi polynomial $P_{k}^{(\alpha,\beta)}$
is solution of the second order differential equation \cite[page
214]{M}
\begin{eqnarray*}  \qquad (1-x^2) y^{''} + [(\alpha-\beta) +(\alpha+\beta+2)x] y^{'} -
k(k+\alpha+\beta+1) y = 0 .\end{eqnarray*} Indeed, by making the
changes $\alpha=a-1, \beta = b-1$ and $k=j+1$, we get
 \begin{eqnarray*}   && \hspace*{-2cm} (x^2-1)\frac{d^2}{d x^2}
 P_{j+1}^{(a-1,b-1)}(x) -[(a-b) +(a+b)x] \frac{d}{d x}
 P_{j+1}^{(a-1,b-1)}(x) + \\ && \hspace*{4cm} + (j+1)(j+a+b)
 P_{j+1}^{(a-1,b-1)}(x) =0.
\end{eqnarray*}
Next, using the fact that  \cite[page 213]{M}
 \begin{eqnarray*}\frac{d}{d x}
   P_{j+1}^{(a-1,b-1)}(x)=
   \frac{j+a+b}{2}P_{j}^{(a,b)}(x),\end{eqnarray*} it follows
    \begin{eqnarray*}
 (x^2-1)\frac{d}{d x}
 P_{j}^{(a,b)} -[(a-b) +(a+b)x]
 P_{j}^{(a,b)}+ 2(j+1)
 P_{j+1}^{(a-1,b-1)} =0.
\end{eqnarray*}
\end{proof}

We conclude this section by giving a sketched proof of Theorem \ref{MainTheorem}.
\begin{proof}[Proof of Theorem \ref{MainTheorem}] We use
the fact that the $L^2$-eigenspace of $\mathfrak{L}_\kappa^{b}$
associated to its first eigenvalue $b$ coincides with the null space
of $\nabla^{*}_{b}$ and is spanned in the Hilbert space
$\mathcal{H}_\kappa$ as follows
$$ Span\set{ (1+\kappa |z|^2)^{-\frac b\kappa} z^n, ~~ n=0,1,2, \cdots } .$$
Hence,  for $b=\nu+  m \kappa$ and according to
(\ref{ClosedEigenfunction}) it turns out to compute $$
\nabla_{\nu+\kappa} \circ \nabla_{\nu+2\kappa}\circ \cdots \circ
\nabla_{\nu + m \kappa}[(1+\kappa |z|^2)^{-(\frac{\nu}{\kappa}+m)}z
^{n}]=: (1+\kappa |z|^2)^{-(\frac\nu\kappa+m)} P_{m,n}^{\nu;
\kappa}(z, \bar z) $$ which is given  by (\ref{RTF3}) in Proposition
\ref{GhanmiProp} and that we can rewrite also as follow
$$
P_{m,n}^{\nu; \kappa}(z, \bar z)
 = (-1)^{m} (m \wedge n)! z^p {\bar z}^q
\mathrm{P}^{(p+q,-2(\frac{\nu}{\kappa}+m)-1)}_{m-q}(1+2\kappa |z|^2)
$$
with $ p=n-m $ and $ q=0$ if  $ m\leq n$ and $ p=0 $  and $q=m-n $
if $ m\geq n$. Next, by applying Proposition \ref{PropLH}, we obtain
the desired result as asserted in Theorem \ref{MainTheorem}.
\end{proof}

\section{ Concluding remarks.}

According to Proposition \ref{GhanmiProp} above, the class of two
variable orthogonal polynomials $ P_{m,n}^{\nu; \kappa}(z, \bar z)$,
suggested by the twisted Laplacians $\mathfrak{L}_\kappa^\nu$ are
obtained in a unified manner taking into account the curvature of
the considered SCRS $M_\kappa$. Indeed the functions $$(1+\kappa
|z|^2)^{-(\frac \nu \kappa +m)} P_{m,n}^{\nu; \kappa}(z, \bar z)$$
are $L^2$-eigenfunctions of $\mathfrak{L}_\kappa^\nu$ with
$E_{\kappa,m}^\nu:= \nu(2m+1)+m(m+1) \kappa$ as corresponding
$L^2$-eigenvalue.  The obtained two variable polynomials satisfy
(\ref{RTF1}),
 where the involved intertwining  invariant
operator $$[(1+ \kappa |z|^2)^{2}\frac{\partial}{\partial z}]^m$$
depends only on the geometry of $M_\kappa$  and not on the magnetic
field. Also, they satisfy the Rodrigues formula (\ref{RTF2}), up to
a given multiplicative constant $C_{\kappa,\nu}^{m,n}$.
  Such polynomials $P_{m,n}^{\nu; \kappa}(z, \bar z)$ are
connected  to the real  Jacobi polynomials
$\mathrm{P}^{(\alpha,\beta)}_{l}(x) $.
Note that the identity (\ref{RTF3}) above for $\kappa\ne 0$ can be rewritten in
the following form
\begin{eqnarray*} P_{m,n}^{\nu; \kappa}(z, \bar z)
=
 (-1)^{m}\left\{
\begin{array}{ll} m! z^{n-m}
\mathrm{P}^{(n-m,-2(\frac{\nu}{\kappa}+m)-1)}_{m}(1+2\kappa z\bar z)
&
\mbox{ if } m\leq n \\
n! {\bar z}^{m-n}
\mathrm{P}^{(m-n,-2(\frac{\nu}{\kappa}+m)-1)}_{n}(1+2\kappa z\bar z)
& \mbox{ if } m\geq n \end{array} \right. \label{RTF4}
  \end{eqnarray*}
Therefore, one sees that the polynomials $P_{m,n}^{\nu; \kappa}(z,
\bar z)$  reduced when $\kappa=-1$ $P_{m,n}^{\nu;
-1}(z, \bar z)$ are exactly the so-called disc polynomials
\cite{Koornw75, Wunsch05}. Here they appear, up to multiplicative functions, as
$L^2$-eigenfunctions of the twisted Laplacian $\mathfrak{L}_{-1}^\nu$ on the hyperbolic disc.  
 While for the limit case $\kappa = 0$, the polynomials  $P_{m,n}^{\nu; 0}(z, \bar z)$ are exactly
  the complex Hermite polynomials \cite{Sh87,IntInt06}
defined by \begin{eqnarray} \label{ComplexeHermite} H_{m,n}^\nu(z,\bar z) := \frac{(-1)^{m+n}}{(2\nu)^n}
e^{2\nu |z|^2}\frac{\partial^{m+n}}{\partial z^{m}\partial\bar
z^{n}}e^{-2\nu |z|^2},\end{eqnarray}
 which form a complete orthogonal system in $L^2(\C;e^{-2\nu |z|^2}dxdy)$.
 The associated functions $e^{-\nu |z|^2}H_{m,n}^\nu(z,\bar z)$ are $L^2$-eigenfunctions of the usual
twisted Laplacian on the Euclidean plane,
$$ \mathfrak{L}_0^\nu= -\{ \frac {\partial^2}{\partial z\partial \bar
z} +\nu (z\frac {\partial}{\partial z}- {\bar
z}\frac{\partial}{\partial{\bar z}}) - \nu^2 |z|^2\},$$ with
$\nu(2m+1)$ as corresponding $L^2$-eigenvalue.

Furthermore, added to the facts that  $\nu(2m+1)=\lim\limits_{\kappa
\longrightarrow 0} E_{\kappa,m}^\nu$ and the operator
$\mathfrak{L}_0^\nu$ appears also as the formal limit of the
unbounded differential operators $\mathfrak{L}_\kappa^\nu$  by
letting $\kappa$ goes to $0$, one gets the limiting transition of $P_{m,n}^{\nu; \kappa}(z, \bar z)$ to the complex
Hermite polynomials $H_{m,n}^\nu(z,\bar z)$:
$$ \lim\limits_{\kappa \longrightarrow 0}P_{m,n}^{\nu;
\kappa}(z, \bar z) =H_{m,n}^\nu(z,\bar z)$$ for every fixed $z\in \C$.
 This can be checked easily from (\ref{RTF2}) or also from (\ref{RTF3})
 using some known useful transformations on special functions.
But below, we give an alternative proof and we will see how this can
be handled using the background related to the factorization method
without knowing explicit expression of $P_{m,n}^{\nu; \kappa}(z,
\bar z)$. Indeed,  by making formal
limit, keeping in mind that $$ \lim\limits_{\kappa \longrightarrow
0}(1+\kappa z\bar z)^{\frac\nu\kappa+m} = e^{\nu z\bar z},
$$ we get
\begin{align*}   \lim\limits_{\kappa \longrightarrow 0}
P_{m,n}^{\nu; \kappa}(z, \bar z) &= \lim\limits_{\kappa
\longrightarrow 0} (1+\kappa z\bar z)^{\frac\nu\kappa+m}
\nabla_{\nu+\kappa} \circ \nabla_{\nu+2\kappa}\circ \cdots \circ
\nabla_{\nu +m\kappa}[(1+\kappa z\bar z)^{-(\frac{\nu}{\kappa}+m)}z
^{n}] \\
&= e^{\nu z\bar z} \underbrace{\nabla_{\nu} \circ \nabla_{\nu}\circ
\cdots \circ \nabla_{\nu}}_{\mbox{m-times}}[e^{-\nu z\bar z}z ^{n}].
\end{align*}
Next, by rewriting $\nabla_{\nu} $ in the following form $
\nabla_{\nu}f = - e^{\nu z\bar z}\frac{\partial }{\partial z}(
e^{-\nu z\bar z}f),$ it follows
\begin{eqnarray*} \lim\limits_{\kappa \longrightarrow 0}
P_{m,n}^{\nu; \kappa}(z, \bar z) = (-1)^m e^{2\nu z\bar z}
\frac{\partial^m}{\partial z^m}[e^{-2\nu z\bar z}z ^{n}].
\end{eqnarray*}
But since $ \frac{\partial^n}{\partial \bar z^n}(e^{-2\nu z\bar z})=
(-2\nu)^n e^{-2\nu z\bar z}z ^{n},$ we conclude easily that
\begin{eqnarray*}  \qquad \lim\limits_{\kappa \longrightarrow 0}
P_{m,n}^{\nu; \kappa}(z, \bar z) = \frac{(-1)^{m+n}}{(2\nu)^n}
e^{2\nu z\bar z} \frac{\partial^{m+n}}{\partial z^m\partial \bar
z^n}(e^{-2\nu z\bar z}) \stackrel{(\ref{ComplexeHermite})}{=}
 H_{m,n}^\nu(z,\bar z) .
\end{eqnarray*}

\quad

 {\bf\it Acknowledgements.} {The author would like to thank all members of CAMS for their hospitality during the full
academic year 2006-2007. He gratefully acknowledges the financial support of the  Arab Regional Fellows
        Program during this period.}

\end{document}